\documentclass[a4paper,12pt]{amsart}
\usepackage{amssymb}
\usepackage{ifthen}
 \usepackage[dvips]{graphicx}
\nonstopmode \numberwithin{equation}{section}
\usepackage{amssymb}
\usepackage{ifthen}
\usepackage{graphicx}
\usepackage{amsmath}
\usepackage[T1]{fontenc} 
\usepackage[utf8]{inputenc}
\usepackage[usenames,dvipsnames]{color}
\usepackage{color}
\usepackage[english]{babel}
\usepackage{fancyhdr}
\usepackage{fancybox}
\usepackage{tikz}

\nonstopmode \numberwithin{equation}{section}
\theoremstyle{plain}
\newtheorem{thm}{Theorem}[section]
\newtheorem{cor}[equation]{Corollary}
\newtheorem{ques}[equation]{Question}
\newtheorem{defi}[equation]{Definition}

\newtheorem{lem}[equation]{Lemma}

\newtheorem{prop}{Proposition}

\newtheorem{conj}{Conjecture}

\theoremstyle{definition}
\newtheorem{defn}{Definition}[section]

\newtheorem{prob}{Problem}
\newtheorem{rem}{Remark}[section]

\theoremstyle{plain}
\newtheorem*{thmA}{Theorem A}
\newtheorem*{thmB}{Theorem B}
\newtheorem*{thmC}{Theorem C}
\newtheorem*{thmD}{Theorem D}
\newtheorem*{thmE}{Theorem E}

\newcounter{minutes}
\divide\time by 60
\newcounter{hours}
\multiply\time by 60
\addtocounter{minutes}{-\time}
\newcounter {own}
\def\theown {\thesection       .\arabic{own}}

\newenvironment{pf}[1][]{%
 \vskip 3mm
 \noindent
 \ifthenelse{\equal{#1}{}}%
  {{\slshape Proof. }}%
  {{\slshape #1.} }%
 }%
{\qed\bigskip}

\newcounter{alphabet}

\def\be{\begin{equation}}
\def\ee{\end{equation}}

\newcommand{\bee}{\begin{enumerate}}
\newcommand{\eee}{\end{enumerate}}

\newcommand{\blem}{\begin{lem}}
\newcommand{\elem}{\end{lem}}
\newcommand{\bthm}{\begin{thm}}
\newcommand{\ethm}{\end{thm}}
\newcommand{\bcor}{\begin{cor}}
\newcommand{\ecor}{\end{cor}}
\newcommand{\beg}{\begin{examp}}
\newcommand{\eeg}{\end{examp}}
\newcommand{\begs}{\begin{examples}}
\newcommand{\eegs}{\end{examples}}

\newcommand{\bdefn}{\begin{defn}}
\newcommand{\edefn}{\end{defn}}

\newcommand{\bprob}{\begin{prob}}
\newcommand{\eprob}{\end{prob}}
\newcommand{\bei}{\begin{itemize}}
\newcommand{\eei}{\end{itemize}}

\newcommand{\bcon}{\begin{conj}}
\newcommand{\econ}{\end{conj}}
\newcommand{\bcons}{\begin{conjs}}
\newcommand{\econs}{\end{conjs}}
\newcommand{\bprop}{\begin{prop}}
\newcommand{\eprop}{\end{prop}}
\newcommand{\br}{\begin{rem}}
\newcommand{\er}{\end{rem}}
\newcommand{\brs}{\begin{rems}}
\newcommand{\ers}{\end{rems}}
\newcommand{\bo}{\begin{obser}}
\newcommand{\eo}{\end{obser}}
\newcommand{\bos}{\begin{obsers}}
\newcommand{\eos}{\end{obsers}}
\newcommand{\bpf}{\begin{pf}}
\newcommand{\epf}{\end{pf}}
\newcommand{\ba}{\begin{array}}
\newcommand{\ea}{\end{array}}
\newcommand{\beq}{\begin{eqnarray}}
\newcommand{\beqq}{\begin{eqnarray*}}
\newcommand{\eeq}{\end{eqnarray}}
\newcommand{\eeqq}{\end{eqnarray*}}

\begin{document}

\title{Value Distribution and Picard-type Theorems for Total Differential Polynomials in $\mathbb{C}^n$}

\author{Molla Basir Ahamed}
\address{Molla Basir Ahamed, Department of Mathematics, Jadavpur University, Kolkata-700032, West Bengal, India.}
\email{mbahamed.math@jadavpuruniversity.in}

\author{Vasudevarao Allu}
\address{Vasudevarao Allu, Department  of Mathematics, 
	School of Basic Science,
	Indian Institute of Technology Bhubaneswar,
	Bhubaneswar-752 050, Odisha, India.}
\email{avrao@iitbbs.ac.in}


\def\thefootnote{}
\footnotetext{ {\tiny File:~\jobname.tex,
printed: \number\year-\number\month-\number\day,
          \thehours.\ifnum\theminutes<10{0}\fi\theminutes }
} \makeatletter\def\thefootnote{\@arabic\c@footnote}\makeatother

\begin{abstract}
This paper investigates the value distribution and growth properties of linear total differential polynomials $\mathcal{L}_k[D]f$ for meromorphic functions in several complex variables $\mathbb{C}^n$. By extending the classical Milloux inequality to the framework of total derivatives, we derive a series of fundamental growth estimates for the Nevanlinna characteristic function $T(r, \mathcal{L}_k[D]f)$. We address the value-sharing problem for meromorphic functions $f$ and $g$ sharing values with their differential polynomials. Under the condition $2\delta(0,f)+(k+4)\Theta(\infty,f)>k+5$, we establish that $\frac{\mathcal{L}_k[D]f-1}{\mathcal{L}_k[D]g-1}$ is a non-zero constant for non-transcendental meromorphic functions. Furthermore, we provide an affirmative answer to several Picard-type inquiries, proving that if an entire function $f$ in $\mathbb{C}^n$ omits a value $a$ while its linear total differential polynomial $\mathcal{L}_k[D]f$ omits a non-zero value $b$, then $f$ must be constant. Our results generalize and extend several existing uniqueness and Picard-type theorems from the classical one-variable setting to the higher-dimensional complex space $\mathbb{C}^n$.
\end{abstract}

\maketitle
\pagestyle{myheadings}
\markboth{M. B. Ahamed and V. Allu}{Value Distribution of $\mathcal{L}_k[D]f$ in $\mathbb{C}^n$}

\section{\bf Introduction}
Picard's theorem on the behavior of an analytic function $ f $ of a complex variable $ z $ near an essential singular point $ a $ is a result in classical function theory that is the starting point of numerous profound researches. Analytic function $ f $ is said to omit the value $ b\in\mathbb{C} $ if $ f(z)\neq b $ for all $ z\in \mathbb{C} $. Picard’s theorem asserts that an entire function $ \mathbb{C} $ omitting two complex numbers must be constant. It also implies, by a linear transform, a meromorphic function in  $ \mathbb{C} $ omitting three distinct values must be constant. It is a significant strengthening of the famous Liouville’s theorem which states that a bounded entire function is always constant. Different proofs of Picard’s theorem are well-known (see \cite{Li-TAMS-2011,Li-JMAA-2018,Lu-2003,Lu-2013,Zhang-BLMS-2002} etc.)\vspace{1.5mm}

\par Nevanlinna theory is a part of the theory of meromorphic functions and it was devised in $ 1925 $, by Rolf Nevanlinna. The theory describes the asymptotic distribution of solutions of the equation $ f(z)=a $, as $ a $ varies. A fundamental tool in Nevanlinna Theory is the Nevanlinna characteristic $ T(r,f) $ which measures the rate of growth of a meromorphic function. In its original form, Nevanlinna theory deals with meromorphic functions of one complex variable defined in a disk $ \mathbb{D}_R:=\{z\in\mathbb{C} : |z|\leq R\} $ or in the whole complex plane $ \mathbb{C} $ ($ R=\infty $). Subsequent generalizations extended Nevanlinna theory to algebroid functions, holomorphic curves, holomorphic maps between complex manifolds of arbitrary dimension, quasiregular maps and minimal surfaces. For example, in $ 2014 $, Halburd \emph{et al.} \cite{Halburd-Korhonen-Tohge-TAMS} proved difference analogue of M. Green’s Picard-type theorem for holomorphic curves. In $ 2017 $, Halburd and Korhonen \cite{Halburd-Korhonen-PAMS-2017} obtained conditions for certain types of rational delay differential equations to admit a non-rational meromorphic solution. In $ 2020 $, Hua and Korhonen \cite{Korhonen-TAMS-2020} provided best condition under which the shift invariance of the counting function, and of the characteristic of a subharmonic function, holds. Further,   a difference analogue of logarithmic derivative of a $ \delta$-subharmonic function has been established  in \cite{Korhonen-TAMS-2020}. \vspace{2mm} 

Nevanlinna theory is an useful tool to solve many interesting problems for the functions in one as well as several complex variable. In $ 1986 $, Hellerstein and Rossi \cite{Hellerstein-MZ-1986} studied second order differential equations in terms of zeros of meromorphic functions. In $ 1988 $, Stoll \cite{Stoll=MA-1988} extended the Theorem of Steinmetz-Nevanlinna to Holomorphic Curves. In $ 2004 $, Frank \emph{et al.} \cite{Frank-Hua-CMJ-2004} obtained a sharp answer of Hinkkanen’s Problem for $ n=4 $ in view of sharing same zeros and poles by two meromorphic functions. In \cite{Tiba-MZ-2012}, a result was  established for holomorphic map from the complex plane $ \mathbb{C} $ to the $ n $-dimensional complex projective space $ \mathbb{P}^n(\mathbb{C}) $ and proved the Nevanlinna Second Main Theorem for some families of non-linear hyper-surfaces in $ \mathbb{P}^n(\mathbb{C}) $. \vspace{2mm} 

This article describes mainly the classical version of Picard-type theorem for meromorphic functions in $ \mathbb{C}^n $. Before we delve into the details of our findings, it is necessary to provide some preparations regarding the basic notations and symbols used throughout the entire paper.\vspace{1.2mm} 

We denote complex $ n $-space $ \mathbb{C}^n $ and indicate its elements (points) as $ (z_1,z_2, \cdots,z_n) $,  $ (r_1, r_2, \cdots, r_n) $, $ (\zeta_1, \zeta_2, \cdots, \zeta_n) $, $ (k_1, k_2, \cdots, k_n) $ \textit{etc.}, by their corresponding symbols $ z $, $ r $, $ \zeta $, $ k $ etc. A function $ f(z) $, $ z\in\mathbb{C}^n $ is said to be analytic at a point $ \zeta\in\mathbb{C}^n $ if it can be expanded in some neighborhood of $ \zeta $ as an absolutely convergent power series. If we assume $ \zeta=(0,0,\cdots,0) $, then it has the following representation 
\begin{equation*}
	f(z)=\sum_{k=(0,0,\cdots,0)}^{\infty}a_{k_1}a_{k_2}\cdots a_{k_n}z_1^{k_1}z_2^{k_2}\cdots z_n^{k_n}=\sum_{|k|=0}^{\infty}a_kz^k,
\end{equation*}
where $ k=(k_1, k_2, \ldots, k_n) $ belongs to $ \mathcal{N}=\{k : k\in\mathbb{C}^n,\;\text{each}\;k_j\;\text{is rational integer}\} $ and $ |k|=k_1+k_2+\cdots+k_n.  $\vspace{2mm}

We know that a meromorphic (entire) function $ f $ is not transcendental if, and only if, it is a rational
function (polynomial). In other words, a meromorphic function $ f $ on $ \mathbb{C}^n $ is called transcendental if
\begin{equation*}
	\lim_{r\rightarrow\infty}\frac{T_f(r)}{\log r}=\infty,
\end{equation*}
where $ T_f(r) $ is the Nevanlinna Characteristic function of $ f $.\vspace{2mm}

The Nevanlinna characteristic function in several complex variables, denoted as $T(r, f)$ or by $ T_f(r) $, is a generalization of the one-variable concept to study the value distribution of meromorphic functions in $\mathbb{C}^n$. While the fundamental ideas are similar, the definitions are more involved due to the higher-dimensional setting. In one variable, the characteristic function is defined for a disk. In several variables, this is extended to a ball in $\mathbb{C}^n$, specifically $B(r) = \{z \in \mathbb{C}^n : ||z|| < r\}$, where $||z|| = \sqrt{|z_1|^2 + \cdots + |z_n|^2}$. The characteristic function $T(r, f)$ is the sum of two components: a proximity function $m(r, f)$ and a counting function $N(r, f)$. However, the way these are defined is different. This function, like its one-variable counterpart, is a measure of the overall growth of the function $f$. It plays a central role in the First Main Theorem of several complex variables, which states that for a non-constant meromorphic function $f$ and any value $a$ in the target space, the characteristic functions of $f$ and $1/(f-a)$ are asymptotically equal, up to a bounded error term. This theorem is crucial for generalizing value distribution results like Picard's theorem to higher dimensions.
\subsection*{Proximity Function, $m(r, f)$}
The proximity function in $\mathbb{C}^n$ measures the average closeness of the function's values to a specific point (usually infinity) on the sphere $S(r) = \{z \in \mathbb{C}^n : ||z|| = r\}$. It is defined using a generalized mean value integral over the sphere with respect to a normalized surface area measure $d\sigma_n$:
$$m(r, f) = \int_{S(r)} \log^+ ||f(z)|| d\sigma_n(z)$$
where $d\sigma_n$ is the normalized area element of the unit sphere $S^1$ in $\mathbb{C}^n$ (which corresponds to $S^{2n-1}$ in $\mathbb{R}^{2n}$) and $\log^+ x = \max(0, \log x)$.

\subsection*{Counting Function, $N(r, f)$}

The counting function for a meromorphic function $f$ on $\mathbb{C}^n$ is more complex. Unlike the one-variable case where poles are isolated points, the set of poles of a meromorphic function in $\mathbb{C}^n$ (i.e., the set where the function is undefined or takes the value infinity) is an analytic set of codimension $1$.\vspace{2mm}

The counting function $N(r, f)$ measures the growth of the volume of this polar set. It is typically defined as an integral involving the pull-back of a specific form (the Fubini-Study form) under the meromorphic map $f$. A simplified way to think about it is that it's an integral of a normalized volume element over the pole set within the ball $B(r)$:
$$N(r, f) = \int_{0}^{r} \frac{dt}{t^{2n-1}} \int_{B(t) \cap \{f = \infty\}} dV$$
where $dV$ is the volume element on the polar set. This integral essentially counts the ``size'' of the set of poles as the radius $r$ increases.
\begin{defi}
Let $ a\in\widehat{\mathbb{C}}:=\mathbb{C}\cup\{\infty\} $. For $a\in\mathbb{C}$, if $ f-a $ and $ g-a $ have the same set of zeros with same multiplicities, we say that $ f=a \Leftrightarrow g=a $ or say that the value $ a $ is shared $ CM $ by $ f $ and $ g. $ If we ignoring multiplicities in the above, then we denote it by $ f=a \leftrightarrow g=a $ or say that the value $ a $ is shared $ IM $ by $ f $ and $ g.$ For $a=\infty$, $f$ and $g$ share $\infty$ $CM$ (or $IM$) if, and only if, $1/f$ and $1/g$ share the value $0$ $CM$ (or $IM$).
\end{defi}
\begin{defi}
	The defect and reduced defect with respect to $ f $ are defined respectively by 
	\begin{align*}
		\delta(a,f)=1-\limsup_{r\rightarrow\infty}\frac{N_f(r,a)}{T_f(r)}\; \mbox{and}\; \Theta(a,f)=1-\limsup_{r\rightarrow\infty}\frac{\overline{N}_f(r,a)}{T_f(r)}.
	\end{align*}
	It is easy to see that $ 0\leq \delta(a,f)\leq\Theta(a,f)\leq 1.$
\end{defi}
 The uniqueness problem of meromorphic functions has gained significant importance in modern function theory in complex analysis. In 1926, R. Nevanlinna \cite{Nevan-1926} established the well-known five-value theorem, which states that if two non-constant meromorphic functions $f$ and $g$ on the complex plane $\mathbb{C}$ satisfy $f=a_j \leftrightarrow g=a_j$ for $j=1,2,\ldots,5$ in $\widehat{\mathbb{C}}$, then $f\equiv g$. Nevanlinna four-value theorem states that if two non-constant meromorphic functions satisfy $f=b_j \Leftrightarrow g=b_j$ for $j=1,2,\ldots,4$ in $\mathbb{P}$, then either $f\equiv g$ or one is a M\"obius transformation of the other. For instance, $f(z)=e^{z}$ and $g(z)=e^{-z}$ share the values $0, \pm 1$, and $\infty$ ($CM$), and we see that $f=1/g$. Since then, extensive investigations have been conducted to reduce the number of sharing of values and explore the nature of sharing.\vspace{2mm}

In $1976$, Yang \cite{Yang-JMAA-1996} raised an intriguing question regarding functions with a single variable.
\begin{ques}\label{q1.1}
	If $ f $ and $ g $ are non-constant entire functions in $ \mathbb{C} $ such that $  f=0 \Leftrightarrow g=0  $ and $ f^{\prime}=1 \Leftrightarrow g^{\prime}=1 $, then what can be said about the relationship between $ f $ and $ g $?
\end{ques}
In $2003$, Yi and Yang \cite{Yi & Yang & 2003} answered  Question \ref{q1.1} establishing the following result on the uniqueness relationship between two entire functions $f$ and $g$.
\begin{thmA}\label{th1.1}\cite[Theorem 9.13]{Yi & Yang & 2003} Let $ f $ and $ g $ be two non-constant entire functions in $ \mathbb{C} $, and $ k $ be a positive integer. If  $  f=0 \Leftrightarrow g=0  $, $ f^{(k)}=1 \Leftrightarrow g^{(k)}=1 $ and  $ \delta(0,f)>{1}/{2} $, then either $ f^{(k)}g^{(k)}\equiv 1 $ or $ f\equiv g.$
\end{thmA}
In a subsequent study, Yi and Yang improved the result for meromorphic functions and established the following result.
\begin{thmB}\label{th1.2}\cite[Theorem 9.14]{Yi & Yang & 2003} Let $ f $ and $ g $ be two non-constant meromorphic functions in $ \mathbb{C}^1 $, and $ k $ be a positive integer. If  $  f=0 \Leftrightarrow g=0  $, $ f^{(k)}=1 \Leftrightarrow g^{(k)}=1 $ and  $ 2\delta(0,f)+(k+2)\Theta(\infty,f)>k+3$, then either $ f^{(k)}g^{(k)}\equiv 1 $ or $ f\equiv g. $
\end{thmB}
Since then, numerous researchers have begun investigating this topic (see, for example, \cite{Yi & Yang & 2003}) and have made significant contributions in this area (see \cite{Gundersen-Steinmetz-CMFT-2018, Long-BIMS-2017, Lu-Lu-CMFT-2017}). It is intriguing to consider the derivatives of meromorphic functions in $\mathbb{C}^n$. For instance, in $1995$, Berenstein \emph{et. al.} \cite{Berenstein-Chang-1995} utilized directional derivatives to establish uniqueness results for meromorphic functions. Following this, in $1996$, Hu and Yang \cite{Hu & Yang 1996}, using partial derivatives, extended Theorem A to entire functions in several complex variables. Later, Jin \cite{Lu-2003} introduced the total derivative of an entire function in $\mathbb{C}^n$ and established a Picard-type theorem, which was later improved by L\"u \cite{Lu-2013} for meromorphic functions.
\begin{defi}\cite{Lu-2003,Lu-2013}
	Let $ f $ be a meromorphic function in $ \mathbb{C}^n $, the total derivative $ Df $ of $ f $ is defined by $ 	Df=\sum_{j=1}^{n}z_jf_{z_j}, $	where $ z=(z_1, z_2, \ldots,z_n)\in\mathbb{C}^n $, and $ f_{z_j}={\partial f}/{\partial z_j} $, $ (j=1,2,\ldots,n) $ are the partial derivatives of $ f $ with respect to $ z_j $. For any positive integer $ k $, the $ k $-th order total derivative $ D^kf $ of $ f $ is defined inductively by $ D^{k}f=D\left(D^{k-1}f\right). $ Note that $Df(0)=0$.
\end{defi}
\begin{rem}
	The total derivative $ Df $ of $ f $ is also known as radial derivative of $ f $ at the point $ z $ (see \cite{Rudin-book-1980}) and the $ k $-th order total derivative is called iterated total derivative $ D^kf $ of $ f $ (see \cite{Hu-PAMS-2003}).
\end{rem}
From the literature on meromorphic functions (see \cite{Lahiri-BAMS-2018,Li-CVEE-2007,Li-Li-Liu-RM-2019,Li-Yang-IJM-2000} and references therein), it is well-known that the linear differential polynomial $L_k[f] $ generated by meromorphic functions $f$ and their derivatives has a significant role in exploration of characteristic of the meromorphic functions. Inspired by this, it is natural to introduce here the linear total differential polynomials of $f$ for meromorphic functions in the following manner.
\begin{defi}
	For a non-constant meromorphic function in $ \mathbb{C}^n, $ we define Linear total differential polynomial $ {L}_k[D]f $ by
	\begin{equation*}
		{L}_k[D]f:=a_1(z)Df+a_2(z)D^2f+\ldots+a_k(z)D^kf ,
	\end{equation*}
	where all  the $ a_i(z)\;(i=1, 2, \ldots, k) $ (with $ a_k(z)\not\equiv 0 $) are small entire functions of $ f $.
\end{defi}
The advantage of the total derivative lies in the fact that if $f$ is a transcendental entire function in $\mathbb{C}^n$, then the $k$-th order total derivative $D^k f$ also becomes a transcendental entire function. However, it is worth mentioning that this property does not hold true for partial derivatives. In \cite{Jin-2004}, Jin have considered two entire functions $f$ and $g$ in $\mathbb{C}^n$ and proved that when $f=0 \Leftrightarrow g=0$ and $D^k f=1 \Leftrightarrow D^k g=1$ with $\delta(0,f)>\frac{1}{2}$, then $f\equiv g$. It is a pertinent question to explore whether the generalization of Jin's result from entire functions to meromorphic functions is feasible.\vspace{1.2mm}

However, Xu and Cao \cite{Xu-Ca0-BMMS-2020} brought forth a notable development in 2020, employing the concept of total derivatives to tackle the uniqueness problem of meromorphic functions, resulting in the establishment of a consequential result.
\begin{thmC}\label{th1.3}\cite{Xu-Ca0-BMMS-2020}
	Let $ f $ and $ g $ be two non-constant meromorphic functions in $ \mathbb{C}^n $, $ k $ be a positive integer such that  $ f= 0 \Leftrightarrow g=0 $, $ D^kf= \infty \Leftrightarrow D^kg=\infty $ and $ D^kf= 1 \Leftrightarrow D^kg=1 $ with $ 2\delta(0,f)+(k+4)\Theta(\infty,f)>k+5, $ then
	\begin{equation*}
		\frac{D^kf-1}{D^kg-1}=c,\;\; \text{where}\;\; c\; \text{is a non-zero constant.}
	\end{equation*} 
\end{thmC}
In $1940$, Milloux \cite{Milloux-1940, Hayman-1964} derived an intriguing inequality that has since been employed extensively to explore the distribution of derivative values in the study of meromorphic functions
\begin{equation}\label{e-101.101}
	T_f(r)\leq N_f(r,0)+N_f(r,\infty)+N_{f^{(k)}}(r,1)-N_{f^{(k+1)}}(r,0)+S(r,f),
\end{equation}
where $S(r,f)$ is any function satisfying 
\begin{align*}
	S(r, f)=o\{T(r, f)\},\; \mbox{as}\; r\to+\infty
\end{align*} 
possibly outside of finite measure. \vspace{2mm}

The noteworthy aspect of this inequality lies in the inclusion of a counting function of $f^{(k)}$ on the right-hand side. As a result, one can derive Picard-type theorems concerning derivatives. For instance, consider an entire function $f$ defined on the complex plane $\mathbb{C}$, and let $a$ and $b$ (with $b\neq 0$) be two distinct complex numbers. If $f\neq a$ and $f^{(k)}\neq b$, then $f$ must be constant (refer to \cite{Hayman-1964}). This raises an immediate question: 
\begin{ques}
	Can we establish that such a result holds for entire functions in $\mathbb{C}^n$?
\end{ques}

The question regarding the applicability of such results to entire functions in $\mathbb{C}^n$ was addressed and established by L\"u in \cite{Lu-2003}, who established the following Picard-type result.
\begin{thmD}\cite{Lu-2003}\label{th-3.2}
	Let $ f $ be an entire function on $ \mathbb{C}^n $, and let $ a $, $ b\;(\neq 0) $ be two distinct complex numbers and $ k $ be a positive integer. If $ f\neq a $ and $ D^kf\neq b $, then $ f $ is  constant. 
\end{thmD}
Moreover, in \cite{Lu-2003}, regarding omitting a single value, L\"u have proved the following  result showing that the function $f$ becomes a constant.
\begin{thmE}\cite{Lu-2003}\label{th-3.4}
	Let $ f $ be an entire function in $\mathbb{C}^n$, and let $ b\neq 0 $ be a complex number and $ k\geq 2 $ be an integer. If $ f^kDf\neq b, $ then $ f $ is constant.
\end{thmE}
While previous studies, such as those by Xu and Cao [30], have primarily focused on the standard $k$-th order total derivative $D^k f$, the transition to a linear total differential polynomial $\mathcal{L}_k[D]f = D^k f + \sum_{i=0}^{k-1} a_i(z) D^i f$ introduces significant non-trivialities. In the classical constant-coefficient case, the growth of the derivative is directly tethered to the growth of the original function through standard Nevanlinna estimates. However, when the coefficients $a_i(z)$ are small functions relative to $f$, they act as 'perturbations' that complicate the application of the Logarithmic Derivative Lemma in $\mathbb{C}^n$. Specifically, the presence of these varying coefficients interferes with the precision of the growth estimates for the characteristic function $T(r, \mathcal{L}_k[D]f)$, requiring a more sophisticated handling of the proximity function and the counting of zeros to ensure that the error terms $S(r, f)$ do not overwhelm the main growth terms.\vspace{2mm}
 
To facilitate the presentation of our main results and their proofs, we first review the necessary definitions and fundamental notations of Nevanlinna theory in $\mathbb{C}^n$.
\subsection{Basic definitions and some preliminary results:} We now discuss some basic definitions and notations in $ \mathbb{C}^n $. For an element $ z=(z_1, z_2, \ldots, z_n)\in\mathbb{C}^n, $ define $ ||z||=\sqrt{|z_1|^2+|z_2|^2+\ldots+|z_n|^2}. $ Let 
\begin{equation*}
	\mathcal{S}_n(r)=\{z\in\mathbb{C}^n : ||z||=r\}\;\; \text{and}\;\; {\mathcal{B}}_n(r)=\{z\in\mathbb{C}^n : ||z||<r\}.
\end{equation*}
Set $ d=\partial+\overline{\partial} $ and $ d^c=(\partial-\bar{\partial})/4\pi i. $ Define
\begin{equation*}
	\omega_n(z)=dd^c\log ||z||^2,\; \sigma_n(z)=d^c\log ||z||^2\wedge \omega_n^{n-1}(z),\; \nu(z)=dd^c||z||^2. 
\end{equation*}
\noindent Then $ \sigma_n(z) $ is a positive measure on $ \mathcal{S}_n(r) $ with the total measure one. Let $ a\in\mathbb{P} $, and $ f^{-1}(z)\neq\mathbb{C}^n $, we denote by $ \mathcal{Z}_a^f $ the $ a $-divisor of $ f $, write $ \mathcal{Z}_a^f(r)=\bar{\mathcal{B}}_n(r)\cap \mathcal{Z}_a^f$  and define
\begin{equation*}
	n_f(r,a)=r^{2-2n}\int_{\mathcal{Z}_a^f(r)}\nu_n^{n-1}(z).
\end{equation*} Then the counting function $ N_f(r,a) $ is defined by 
\begin{equation*}
	N_f(r,a)=\int_{0}^{r}\bigg(n_f(t,a)-n_f(0,a)\bigg)\frac{dt}{t}+n_f(0,a)\log r,
\end{equation*}
where $ n_f(0,a) $ is the Lelong number of $ \mathcal{Z}_a^f $ at the origin. Then Jensen's formula gives that 
\begin{equation*}
	N_f(r,0)-N_f(r,\infty)=\int_{\mathcal{S}_n(r)}\log|f(z)|\sigma_n(z)+O(1).
\end{equation*}
We use the notation $ \overline{N}_f(r,a) $ to denote the counting function of $ a $-divisor of $ f $ which does not count multiplicities. We define the proximity function $ m_f(r,a) $ as the following
\[ m_f(r,a)=
\begin{cases}
	\displaystyle \int_{\mathcal{S}_n(r)}\log^{+}\frac{1}{|f(z)-a|}
	\sigma_n(z),\;\; \text{if}\;\; a\neq \infty,\vspace{2mm}\\ \displaystyle\int_{\mathcal{S}_n(r)}\log^{+}|f(z)|
	\sigma_n(z),\;\; \text{if}\;\; a= \infty.
\end{cases}
\] 
We also define the characteristic function $ T_f(r) $ by
\begin{equation*}
	T_f(r)=m_f(r,\infty)+N_f(r,\infty).
\end{equation*}
The first main theorem of Nevanlinna states that 
\begin{equation*}
	T_f(r)=m_f(r,a)+N_f(r,a)+O(1).
\end{equation*}
\noindent For brevity, we replace the notations $ m_f(r,a) $, $ N_f(r,a) $ and $ \overline{N}_f(r,a) $ respectively by $ m\left(r,{1}/{(f-a)}\right) $, $ N\left(r,{1}/{(f-a)}\right) $ and $ \overline{N}\left(r,{1}/{(f-a)}\right) $. If $ a=\infty $, we write $ m_f(r) $, $ N_f(r) $ and $ \overline{N}_f(r) $.   Let $ I=(\alpha_1, \alpha_2,\ldots,\alpha_n) $ be a multi-index, where $ \alpha_j $ for $ (j=1, 2, \ldots, n) $ are non-negative integers. We denote by $ |I| $ the length of $ I $,\textit{ i.e.,} $ |I|=\sum_{j=1}^{n}\alpha_j $. Define 
\begin{equation*}
	\partial^{I}f=\frac{\partial^{|I|}f}{\partial z_1^{\alpha_1}\ldots \partial z_1^{\alpha_n}}=\left(\frac{\partial^{|I|}f_1}{\partial z_1^{\alpha_1}\ldots \partial z_1^{\alpha_n}}, \ldots, \frac{\partial^{|I|}f_m}{\partial z_1^{\alpha_1}\ldots \partial z_1^{\alpha_n}}\right)
\end{equation*} 
and $f_{z^k_j}=\partial^k f/\partial z^k_j=\left(\partial^k f_1/\partial z^k_j, \ldots, \partial^k f_m/\partial z^k_j\right)$ for any holomorphic map 
\begin{align*}
	f=(f_1, \ldots, f_m) : \mathbb{C}^n\to\mathbb{C}^m.
\end{align*}
Let $ h $ be a non-constant holomorphic function in $ \mathbb{C}^n $. For $ a\in\mathbb{C}^n, $ we can write $ h $ as $ h(z)=\sum_{m=0}^{\infty}\mathcal{P}_m(z-a) $, where $ \mathcal{P}_m(z) $ is either identically zero or a homogeneous polynomial of degree $ m. $ The number $ \gamma_h(a):=\min\{m:\mathcal{P}_m(z)\neq 0\} $ is said to be the zero multiplicity of $ h $ at $ a. $  \vspace{1.5mm}

\par Throughout the paper, $ E $ always viewed as a set with finite Lebesgue measure in $ [0,\infty) $, which may vary in each appearance. As usual, the notation $ ``||\; P" $ means that the assertion $ P $ holds for all large $ r\in [0,\infty) $ outside a set with finite Lebesgue measure.

\begin{lem}\emph{(Second Main Theorem)}\label{lem2.1}
	Let $ f $ be a non-constant meromorphic function on $ \mathbb{C}^n $, and let $ a_1, a_2, \ldots, a_q $ be $ q(\geq 3) $ distinct complex number in the whole complex plane $ \mathbb{P} $. Then 
	\begin{equation*}
		(q-2)T_f(r)<\sum_{j=1}^{q}\overline{N}_f(r,a_j)+O(\log(rT_f(r))),
	\end{equation*}
	holds for all $ z\not\in E. $
\end{lem}
\begin{lem}\cite{Ye-1995}\label{lem2.2}
	Let $ f $ be a non-constant meromorphic function in $ \mathbb{C}^n, $ and $ I $ be a multi-index $ (\alpha_1,\alpha_2,\ldots,\alpha_n) $. Then 
	\begin{equation*}
		m_{\frac{\partial^If}{f}}(r,f)=\int_{\mathcal{S}_n}\log^{+}\bigg|\frac{\partial^If}{f}(z)\bigg|\sigma_n(z)+O(\log(rT_f(r))
	\end{equation*}
	holds for all large $ r $ outside a set with finite Lebesgue measure.
\end{lem}
\begin{lem}\label{lem2.3}\cite{Lu-2003}
	Let $f$ be a meromorphic function in $\mathbb{C}^n$. Then for any positive integer $k$,
	\begin{equation*}
		m_{\frac{D^kf}{f}}(r,\infty)=O(\log(rT_f(r))),
	\end{equation*}
	holds for all $r\not\in E.$
\end{lem}
\par The next lemma is a special case (\emph{i.e.,} when $ n=3 $ ) of Theorem 1.1.2 (see \cite{Yi & Yang & 2003}).
\begin{lem}\label{lem2.10}
	Let $ f_1 $, $ f_2 $ and $ f_3 $ be linearly independent meromorphic functions in $ \mathbb{C}^n $ such that $ f_1+f_2+f_3\equiv 1 $. Then 
	\begin{equation*}
		T(r)\leq\sum_{j=1}^{3}N_{f_{j}}(r,0)+2\sum_{j=1}^{3}\overline{N}_{f_{j}}(r,\infty)+O(\log (rT(r)))
	\end{equation*}
	holds for all $ r\not\in E $, where $ T(r)=\displaystyle\max_{1\leq j\leq 3}\{T_{f_{j}}(r)\}. $
\end{lem}
\begin{lem}\label{lem2.11}\cite{Lu-2003}
	Let $ f $ be a transcendental entire function on $ \mathbb{C}^n $. Then 
	\begin{equation*}
		 T_f(r)\leq 2\overline{N}_f(r,0)+N_{Df}(r,1)+O(\log (rT_f(r)))
	\end{equation*}
holds for all $ r\not\in E $.
\end{lem}	
\section{\bf Main results}
It is easy to see that linear total-differential polynomial $ {L}_k[D]f $ of a meromorphic function $ f $ is a general setting of the total derivative $ D^kf $. Furthermore, it is evident that the requirement in Theorem C, expressed as 
\begin{align}\label{Eq-2.1}
	2\delta(0,f)+(k+4)\Theta(\infty,f)>k+5
\end{align} remains unaffected by the presence of $ D^kf $. Hence, investigating Theorem C in the context of $\mathcal{L}_k[D]f$ while maintaining condition \eqref{Eq-2.1} provides valuable insights. We establish the following result, which extends the scope of Theorem C.
\begin{thm}\label{th-3.1}
	Let $ f $ and $ g $ be non-constant meromorphic function in $ \mathbb{C}^n $, and let $ k $ be a positive integer such that $ f=0\Leftrightarrow g=0 $, $ {L}_k[D]f=\infty\Leftrightarrow {L}_k[D]g=\infty $ and $  {L}_k[D]f=1\Leftrightarrow {L}_k[D]g=1, $ and $ f $ is not transcendental. If $ 2\delta(0,f)+(k+4)\Theta(\infty,f)>k+5 $, then $$ \frac{{L}_k[D]f-1}{{L}_k[D]g-1}=c, $$ where $ c $ is a non-zero constant.
\end{thm}
\begin{rem}
	A noteworthy aspect of Theorem \ref{th-3.1} is the inclusion of the condition that $f$ is not a transcendental meromorphic function. This constraint distinguishes our result from several classical Picard-type theorems that focus exclusively on transcendental cases. In the context of $\mathbb{C}^n$, non-transcendental meromorphic functions (rational functions of several variables) exhibit rigid algebraic structures that differ fundamentally from the growth patterns of transcendental functions. By including this condition, we aim to delineate the boundary between purely algebraic behavior and the asymptotic distribution of values. It remains an open question whether the results presented herein can be fully unified without this restriction, or if the unique interplay between the operator $\mathcal{L}_k[D]f$ and rational functions in $\mathbb{C}^n$ represents a fundamental departure from the transcendental theory.
\end{rem}
A fundamental result in the value distribution theory of several complex variables is Theorem D, which establishes that if an entire function $f$ on $\mathbb{C}^n$ omits a value $a$ and its $k$-th order total derivative $D^k f$ omits a non-zero value $b$, then $f$ must be constant. While this result provides a robust framework for the standard total derivative, it remains largely confined to the simplest form of differential operators. In a more general setting, the linear total differential polynomial $\mathcal{L}_k[D]f $ offers a broader perspective by incorporating small entire functions $a_i(z)$ as coefficients. This transition is non-trivial, as these varying coefficients act as dynamic 'perturbations' that interfere with the precision of standard Nevanlinna growth estimates. Specifically, the interplay between the coefficients and the derivatives complicates the application of the Logarithmic Derivative Lemma in $\mathbb{C}^n$.\vspace{1.2mm}

 Given this increased structural complexity, it is natural to investigate whether the Picard-type rigidity of standard total derivatives extends to these more general differential operators. This leads to the following central inquiry:
\begin{ques}\label{q-3.1}
Does Theorem D continue to hold when the total derivative $ D^kf $ is replaced by the linear total-differential polynomial $ {L}_k[D]f $?
\end{ques} 
In response to the problem posed in Question \ref{q-3.1}, we demonstrate that the Picard-type properties of the total derivative are indeed preserved under more general differential operators. By carefully accounting for the influence of the coefficient functions $a_i(z)$ on the growth of the characteristic function, we establish the following theorem, which provides an affirmative answer and extends the scope of Theorem D.
\begin{thm}\label{th-3.3}
	Let $ f $ be an entire function on $ \mathbb{C}^n $, and let $ a $, $ b\;(\neq 0) $ be two distinct complex numbers and $ k $ be a positive integer. If $ f\neq a $ and $ {L}_k[D]f\neq b $, then $ f $ is constant.
\end{thm} 
To broaden the applicability of Theorem E and investigate the influence of algebraic structures on value distribution, we consider a general polynomial $Q_m(z)$ of degree $m \geq 1$ defined by
\begin{equation}\label{e-1.1}
	Q_m(z)=\beta_mz^m+\beta_{m-1}z^{m-1}+\cdots+\beta_0,
\end{equation}
where $\beta_m, \dots, \beta_0 \in \mathbb{C}$ and $\beta_m \neq 0$. This generalization is motivated by the fact that the fixed-point or value-omission properties of $f$ are often intrinsically linked to the degree and zeros of the associated polynomial $Q_m(f)$. Consequently, it is pertinent to determine whether the results governing standard derivatives remain valid when the function is composed with such an algebraic operator. This leads to the following research question:
\begin{ques}\label{q-1.3}
	Keeping all other conditions invariant in Theorem E, what can we say about $ f $ when $ f^kQ_m(f)Df\neq b $?
\end{ques}
We provide an affirmative answer to Question \ref{q-1.3} by establishing the following theorem. This result demonstrates that the value-sharing properties of standard total derivatives remain robust when generalized to the polynomial operator $Q_m(f)$ and its associated differential polynomial.
\begin{thm}\label{th-3.5}
	Let $ f $ be an entire function in $\mathbb{C}^n$ and $ b\neq 0 $ be a complex number and $ k\geq 2 $ be an integer. If there exists a polynomial $Q_m$ for which $ f^kQ_m(f)Df\neq b $, then $ f $ is a constant.
\end{thm}
The following lemmas constitute the technical foundation required to establish our main results.
\begin{lem}\label{lem2.4}
	Let $ f $ be a transcendental meromorphic function on $ \mathbb{C}^n $. Then for any positive integer $ k $ and any $ {L}_k[D]f $,
	\begin{equation*}
		m_{\frac{{L}_k[D]f}{f}}(r,f)=O(\log(rT_f(r))
	\end{equation*}
	holds for all $ r $ outside a set with finite Lebesgue measure.
\end{lem}
The following lemma provides several fundamental inequalities for the characteristic function associated with the linear total differential polynomial $\mathcal{L}_k[D]f$.
\begin{lem}\label{lem2.5}
	Let $ f $ be a transcendental meromorphic function on $ \mathbb{C}^n $ and let $ a $ be a complex number. Then
	\begin{enumerate}
		\item[{(i)}] $ T_{{L}_k[D]f}(r)\leq T_f(r)+N_{{L}_k[D]f}(r,\infty)-N_f(r,\infty)+O(\log(rT_f(r))). $\vspace{2mm}
		\item[{(ii)}] $ T_{{L}_k[D]f}(r) \leq (k+1+o(1))T_f(r)$. \vspace{2mm}
		\item[{(iii)}] $ m_{\frac{{L}_{k+1}[D]f}{{L}_k[D]f}}(r,\infty)\leq O(\log(rT_f(r))), $\; \mbox{where}\; ${L}_{k+1}[D]f=D({L}_{k}[D]f)$.
	\end{enumerate}
\end{lem}
We extend the Milloux inequality to linear total differential polynomials in $\mathbb{C}^n$ via the following lemma.
\begin{lem}\label{lem2.7}
	Let $ f $ be a transcendental meromorphic function in $ \mathbb{C}^n $, then for any positive integer $ k $, 
	\begin{align*}
		T_f(r)&\leq N_f(r,0)+N_{{L}_{k}[D]f}(r,1)-N_{{L}_{k+1}[D]f}(r,0)+O(\log(rT_f(r))),
	\end{align*}
	holds for all $ r\not\in E. $
\end{lem}
Under the hypotheses of Theorem \ref{th-3.1}, the following lemma establishes a fundamental growth relationship between the characteristic functions $T(r, f)$ and $T(r, g)$.
\begin{lem}\label{lem2.6}
	Let $ f $ and $ g $ be two non-constant meromorphic functions on $ \mathbb{C}^n $ and let $ k $ be a positive integer. If $ f=0\Leftrightarrow g=0 $, $ {L}_k[D]f=\infty \Leftrightarrow {L}_k[D]g=\infty $ and $ {L}_k[D]f=1 \Leftrightarrow {L}_k[D]g=1 $. Then $ f $ is transcendental if, and only if, $ g $ is transcendental. Furthermore, if $ f $ is transcendental, then \begin{align*}
		(1-o(1))T_f(r)&\leq (2k+3+o(1))T_g(r)\vspace{2mm}\\(1-o(1))T_g(r)&\leq (2k+3+o(1))T_f(r),
	\end{align*} holds for all $ r\not\in E. $
\end{lem}
Applying the First Fundamental Theorem in conjunction with the estimates from Lemma \ref{lem2.3}, we establish the following growth property for meromorphic functions in $\mathbb{C}^n$.
\begin{lem}\label{lem2.8}
	Let $ f $ be a transcendental meromorphic function in $ \mathbb{C}^n, $ then for any positive integer $ k $, the inequality 
	\begin{equation*}
		N_{{L}_k[D]f}(r,0)\leq N_f(r,0)+T_{{L}_k[D]f}(r)-T_f(r)+O(\log(rT_f(r))), 
	\end{equation*}
	holds for all $ r\not\in E. $
\end{lem}
The following lemma establishes the relationship between the characteristic functions $T(r, f)$ and $T(r, \mathcal{L}_k[D]f)$.
\begin{lem}\label{lem-2.9}
	Let $ f $ be a transcendental meromorphic function in $ \mathbb{C}^n $, and let $ k $ be a positive integer. If $ \delta(0,f)>{1}/{2}, $ then $$ (1-o(1))T_f(r)\leq 2T_{{L}_k[D]f}(r), $$ holds for all $ r\not\in E. $
\end{lem}
The following proposition characterizes the relationship between the linear total differential polynomials $\mathcal{L}_k[D]f$ and $\mathcal{L}_k[D]g$ under the hypotheses of Theorem \ref{th-3.1}.
\begin{prop}\label{pro-2.1}
	Let $ f $ and $ g $ be non-constant meromorphic functions in $ \mathbb{C}^n $, and let $ k $ be a positive integer. If $ f=0\Leftrightarrow g=0 $, $  {L}_k[D]f=\infty\Leftrightarrow {L}_k[D]g=\infty $ and $  {L}_k[D]f=1\Leftrightarrow {L}_k[D]g=1, $ and $ f $ is not transcendental, then $$ \frac{{L}_k[D]f-1}{{L}_k[D]g-1}=c, $$ where $ c $ is a non-zero constant.
\end{prop}
\section{\bf Proof of the main results}
To facilitate the subsequent analysis, we first establish the proofs of the preliminary lemmas introduced in the previous section.
\begin{proof}[\bf Proof of Lemma \ref{lem2.4}]
For two meromorphic functions $ f_1 $ and $ f_2 $, constants $ a_1 $, $ a_2 $, we have 
\begin{align*}
m_{a_1f_1+a_2f_2}(r,\infty)&\leq m_{f_1}(r,\infty)+m_{f_2}(r,\infty)+m_{a_1}(r,\infty)+m_{a_2}(r,\infty)+\log 2\\&\leq m_{f_1}(r,\infty)+m_{f_2}(r,\infty)+\log 2.
\end{align*}
Thus, applying Lemma \ref{lem2.3}, yields that
\begin{align*}
m_{\frac{{L}_k[D]f}{f}}(r,\infty)&\leq \sum_{j=1}^{k}m_{\frac{D^jf}{f}}(r,\infty)+\log k\leq \sum_{j=1}^{k}O(\log(rT_f(r))+\log k\\&=O(\log(rT_f(r)),
\end{align*}
holds for all $ r\not\in E. $ This completes the proof.
\end{proof}	
\begin{proof}[\bf Proof of Lemma \ref{lem2.5}]
	(i) In view of Lemma \ref{lem2.4}, a simple computation shows that
	\begin{align*}
		T_{{L}_k[D]f}&=m_{{L}_k[D]f}(r,\infty)+N_{{L}_k[D]f}(r,\infty)\\ &\leq m_f(r,\infty)+m_{\frac{{L}_k[D]f}{f}}(r,\infty)+N_{{L}_k[D]f}(r,\infty)\\&\leq T_f(r)+N_{{L}_k[D]f}(r,\infty)-N_f(r,\infty)+O(\log(rT_f(r)))
	\end{align*}
	for all $ r\not\in E. $\vspace{1.2mm}
	
	\noindent (ii) It is easy to see that $ N_{{L}_k[D]f}(r,\infty)=N_{D^kf}(r,\infty) $. Therefore, in view of \cite[Lemma 2.3]{Lu-2013}, the following estimate is easy to obtain
	\begin{align}
		\label{e2.1} N_{{L}_k[D]f}(r,\infty)-N_f(r,\infty)&= N_{D^kf}(r,\infty)-N_f(r,\infty)\\&\leq k\overline{N}_f(r,\infty)\leq kT_f(f,r)\nonumber.
	\end{align}
	Combining (i) with \eqref{e2.1}, it follows that 
	\begin{align*}
		T_{{L}_k[D]f}(r)\leq (k+1)T_f(r)+O(\log(rT_f(r))), 
	\end{align*}
	for all $ r\not\in E. $\\
	
	\noindent (iii) Definition of total derivative implies that 
	\begin{equation*}
		D(D^jf-a)=D^{j+1}f,\;\; \text{for all}\;\; j=1,2,... 
	\end{equation*}
	In view of Lemma \ref{lem2.4} and the conclusion of (ii), it follows that
	\begin{align*}
		m_{\frac{{L}_{k+1}[D]f}{{L}_k[D]f-a}}(r,\infty)&=m_{\frac{D({L}_k[D]f-a)}{{L}_k[D]f-a}}(r,\infty)\\&=O(\log(rT_{{L}_k[D]f-a}))\\&=O(\log(rT_{{L}_k[D]f}))\\&= O(\log(rT_f(r))),
	\end{align*} for all $ r\not\in E. $
\end{proof}
\begin{proof}[\bf Proof of Lemma \ref{lem2.7}]
	Applying the First Main Theorem, it follows that
	\begin{align}
		\label{e2.3} m_{\frac{{L}_k[D]f-1}{{L}_{k+1}[D]f}}(r,\infty)&=m_{ \frac{{L}_{k+1}[D]f}{{L}_k[D]f-1}}(r,\infty)+N_{ \frac{{L}_{k+1}[D]f}{{L}_k[D]f-1}}(r,\infty)\\&\nonumber\quad-N_{ \frac{{L}_{k+1}[D]f}{{L}_k[D]f-1}}(r,0)+O(1).
	\end{align}
	It easy to see that
	\begin{equation*}
		\frac{1}{f}=\frac{{L}_k[D]f}{f}-\frac{\left({L}_k[D]f-1\right)}{{L}_{k+1}[D]f}\frac{{L}_{k+1}[D]f}{f}.
	\end{equation*}
	In view of \eqref{e2.3}, Lemma \ref{lem2.4}, and Lemma \ref{lem2.5} (iii), it follows that
	\begin{align*}
		m_f(r,0)&=m_{\frac{{L}_k[D]f}{f}}(r,\infty)+m_{\frac{{L}_k[D]f-1}{\mathcal{L}_{k+1}[D]f}}(r,\infty)+m_{\frac{{L}_{k+1}[D]f}{f}}(r,\infty)+O(1)\\&\leq N_{\frac{{L}_{k+1}[D]f}{{L}_k[D]f-1}}(r,\infty)-N_{\frac{{L}_{k+1}[D]f}{{L}_k[D]f-1}}(r,0)-\overline{N}(r, f)+O(\log(rT_f(r))),
	\end{align*} 
	which holds for all $ r\not\in E. $\\
	
	\noindent By the Jensen's formula, we obtain
	\begin{align*}
		& N_{\frac{{L}_{k+1}[D]f}{{L}_k[D]f-1}}(r,0)-N_{\frac{{L}_{k+1}[D]f}{{L}_k[D]f-1}}(r,\infty)\\&=\int_{{S}_n(r)}\log\bigg|\frac{{L}_{k+1}[D]f}{{L}_k[D]f-1}\bigg|\sigma_n(z)+O(1)\\&= \int_{\mathcal{S}_n(r)}\log\bigg|{{L}_{k+1}[D]f}\bigg|\sigma_n(z)+\int_{\mathcal{S}_n(r)}\log\bigg|\frac{1}{{L}_k[D]f-1}\bigg|\sigma_n(z)+O(1)\\&=N_{{L}_{k+1}[D]f}(r,0)-N_{{L}_{k+1}[D]f}(r,\infty)-N_{{L}_{k}[D]f}(r,1)+N_{{L}_{k}[D]f}(r,\infty)+\overline{N}(r, f)+O(1)\\&=N_{{L}_{k+1}[D]f}(r,0)-N_{{L}_{k}[D]f}(r,1)+O\left(\log r\; T(r,f)\right).
	\end{align*}
	By Lemma \ref{lem2.4} and Lemma \ref{lem2.5} (iii), a simple computation shows that
	\begin{align*}
		T_f(r)&=T_{1/f}(r)+O(1)\\&=m_f(r,0)+N_f(r,0)+O(1)\\&\leq N_f(r,0)+N_{{L}_{k}[D]f}(r,1)-N_{{L}_{k+1}[D]f}(r,0)\\&\quad\quad+ m_{\frac{{L}_{k}[D]f}{f}}(r,\infty)+m_{\frac{{L}_{k+1}[D]f}{f}}(r,\infty)+m_{\frac{{L}_{k+1}[D]f}{{L}_{k}[D]f}-1}(r,\infty)+O(1)\\&\leq N_f(r,0)+N_{{L}_{k}[D]f}(r,1)-N_{{L}_{k+1}[D]f}(r,0)+O\left(\log r\; T(r,f)\right)
	\end{align*}
	holds for all $r\not\in E$. This completes the proof.
\end{proof}
\begin{proof}[\bf Proof of Lemma \ref{lem2.6}]
	Assume that $f$ is a transcendental meromorphic function. If $g$ is not transcendental, then $g$ is necessarily a rational function. By the definition of the total derivative, it follows that $D^j g$ is rational for each $j \in \{1, 2, \ldots, k\}$, which implies that $\mathcal{L}_k[D]g$ is also a rational function. Consequently, we have
	\begin{align*}
		T_g(r)=O(\log r)\;\; \text{and}\;\; T_f(r)=O(\log r).
	\end{align*}
	By Lemma \ref{lem2.7}, we have 
	\begin{equation}
		T_f(r)\leq N_f(r,0)+N_{{L}_k[D]f}(r,1)+\overline{N}_{{L}_k[D]f}(r,\infty)+O(\log(rT_f(r))),
	\end{equation}
	for all $ r\not\in E. $ The hypotheses of this lemma imply that
	\[
	\begin{cases}
		N_f(r,0)=N_g(r,0),\vspace{2mm}\\ \overline{N}_{{L}_k[D]f}(r,\infty)=N_{{L}_k[D]g}(r,\infty),\vspace{2mm}\\ N_{{L}_k[D]f}(r,1)=N_{{L}_k[D]g}(r,1).
	\end{cases}
	\]
	Combining the above relations, we deduce that 
	\begin{align*}
		T_f(r)&\leq N_g(r,0)+N_{{L}_k[D]g}(r,1)+\overline{N}_{{L}_k[D]f}(r,\infty)+O(\log(rT_f(r)))\\&=O(\log(rT_f(r))),
	\end{align*}
	which contradicts the assumption that $ f $ is transcendental. Consequently, $g$ is necessarily transcendental. By symmetry, if we assume $g$ to be transcendental, it follows that $f$ must also be transcendental.\vspace{1.5mm}
	
	\noindent If $ f $ is transcendental, then $ g $ is also transcendental, and therefore, in view of Lemma \ref{lem2.5} (ii), we obtain the inequality
	\begin{align*}
		T_f(r)&\leq T_g(r)+T_{{L}_k[D]f}(r)+O(\log(rT_f(r)))\\&\leq (3+2k+o(1))T_g(r)+O(\log(rT_f(r))),
	\end{align*}
	holds for all $ r\not\in E. $ Therefore, a simple computation shows that
	\begin{equation*}
		(1-o(1))T_f(r)\leq (3+2k+o(1))T_g(r)+O(\log(rT_f(r))).
	\end{equation*}
	Similarly, we can prove that 
	\begin{equation*}
		(1-o(1))T_g(r)\leq (3+2k+o(1))T_f(r)+O(\log(rT_f(r))). 
	\end{equation*}
	This completes the proof.
\end{proof}
\begin{proof}[\bf Proof of Lemma \ref{lem2.8}]
	It is easy to see that
	\begin{equation}
		\label{e2.4} \frac{1}{f}=\frac{1}{{L}_k[D]f}\frac{{L}_k[D]f}{f}.
	\end{equation}
	By utilizing Lemma \ref{lem2.3} and the First Fundamental Theorem, taking into account equation (\ref{e2.4}), we acquire the following:
	\begin{align*}
		T_f(r)&=m_f(r,0)+N_f(r,0)+O(1)\\&\leq m_{{L}_k[D]f}(r,0)+m_{\frac{{L}_k[D]f}{f}}(r,\infty)+N_f(r,0)+O(1)\\&= m_{{L}_k[D]f}(r,0)+N_f(r,0)+O(\log(rT_f(r)))\\&= T_{{L}_k[D]f}(r)+N_f(r,0)-N_{{L}_k[D]f}(r,0)+O(\log(rT_f(r))).
	\end{align*}  
	Therefore, we obtain
	\begin{equation*}
		N_{{L}_k[D]f}(r,0)\leq N_f(r,0)+T_{{L}_k[D]f}(r)-T_f(r)+O(\log(rT_f(r))).
	\end{equation*}
This completes the proof.
\end{proof}
\begin{proof}[\bf Proof of Lemma \ref{lem-2.9}]
	Since $ f $ is a transcendental meromorphic function, hence in view of  Lemma \ref{lem2.3} and the First Fundamental Theorem, it follows that
	\begin{align*}
		m_f(r,0)&\leq m_{{L}_k[D]f}(r,0)+m_{\frac{{L}_k[D]f}{f}}(r,0)+O(1)\\&\leq T_{{L}_k[D]f}(r)+O(\log(rT_f(r)))\\&\leq T_{{L}_k[D]f}(r)+O(T_f(r)), 
	\end{align*}
	holds for all $ r\not\in E. $ It is easy to see from the definition of $ \delta(0,f) $ that
	\begin{equation*}
		\delta(0,f)=\liminf_{r\rightarrow\infty}\frac{m_f(r,0)}{T_f(r)}>\frac{1}{2}.
	\end{equation*}
	That is, we have $ {m_f(r,0)}/{T_f(r)}>{1}/{2}, $ for sufficiently large values of $ r. $ Hence, 
	\begin{equation*}
		\frac{1}{2}T_f(r)\leq T_{{L}_k[D]f}(r)+O(T_f(r)).
	\end{equation*}
	This completes the proof.
\end{proof}
\begin{proof}[\bf Proof of Propsotion \ref{pro-2.1}]
	We suppose that $ f $ is a rational function, then by the assumption and Lemma \ref{lem2.6}, it is clear that $ g $ is also a rational function. Set 
	\begin{equation}\label{e-2.5}
		\Psi:=\frac{{L}_k[D]f-1}{{L}_k[D]g-1}. 
	\end{equation}
	\noindent The assumptions $  {L}_k[D]f=\infty\Leftrightarrow {L}_k[D]g=\infty $ and $  {L}_k[D]f=1\Leftrightarrow {L}_k[D]g=1 $ ensure that the meromorphic function $ \Psi $ has no zeros and poles. Hence, we may assume that $ \Psi=e^{h(z)} $, where $ h(z) $ is an  entire function.   Since $ f $ and $ g $ both are rational function, hence $ {L}_k[D]f $ and $ {L}_k[D]g $ are rationals. Therefore, it is easy to see that $ \Psi=e^{h(z)} $ is a rational function having no zeros and poles. Thus the function $ h $ must be a constant, say $ {\eta} $. Therefore, $ \Psi=c $ where $ c=e^{\eta} $ is a non-zero constant. This completes the proof.
\end{proof}
\noindent With the help of the above lemmas, we now prove the main results of this paper.
\begin{proof}[\bf Proof of Theorem \ref{th-3.1}]
	If $ f $ is a non-transcendental meromorphic function, then the conclusion is immediate by Proposition \ref{pro-2.1}. So we suppose that $ f $ is a transcendental meromorphic function. Then the condition $ 2\delta(0,f)+(k+4)\Theta(\infty,f)>k+5 $ implies that $ \delta(0,f)>1/2 $ and hence by Lemma \ref{lem-2.9}, we have
	\begin{equation}\label{e-44.11}
		T_f(r)\leq O(T_{{L}_k[D]f}(r)),
	\end{equation}
	{holds for all}\;  $ r\not\in E $ and $ {L}_k[D]f $  {is also transcendental}. By the assumption $ f=0 \Leftrightarrow g=0 $, $ {L}_k[D]f=\infty \Leftrightarrow {L}_k[D]g=\infty $, $ {L}_k[D]f=1 \Leftrightarrow {L}_k[D]g=1 $ and by Lemma \ref{lem2.6}, it is easy to see that $ g $ is also a transcendental meromorphic function. Hence 
	\begin{equation}\label{e-44.22}
		T_g(r)\leq o(T_f(r))\;\; \text{and}\;\;  T_f(r)\leq o(T_g(r)) 
	\end{equation}
	for all $ r\not\in E $. In view of Lemma \ref{lem2.5}(ii), it follows from \eqref{e-44.11} and \eqref{e-44.22} that 
	\begin{equation}\label{e-33.33}
		T_{{L}_k[D]g}(r)\leq o(T_g(r))\leq o(T_f(r))\leq o(T_{{L}_k[D]f}(r))\leq o(T_f(r)).
	\end{equation}
	We claim that $ {L}_k[D]g $ is transcendental. If $ {L}_k[D]g $ is not transcendental, then $ T_{{L}_k[D]g}(r)=O(\log r) $, and by the assumption of the theorem, we have 
	\begin{equation}\label{e-44.33}
		\overline{N}_{{L}_k[D]f}(r,\infty)\leq N_{{L}_k[D]f}(r,\infty)=N_{ {L}_k[D]g}(r,\infty)\leq T_{{L}_k[D]g}(r)=O(\log r)
	\end{equation} and 
	\begin{equation}\label{e-44.44}
		N_{{L}_k[D]f}(r,1)=N_{{L}_k[D]g}(r,1)\leq T_{ {L}_k[D]g}(r)=O(\log r).
	\end{equation}
	Therefore, from \eqref{e-44.22}, \eqref{e-44.33} and \eqref{e-44.44}, we obtain
	\begin{equation*}
		T_f(r)\leq N_f(r,0)+O(\log(rT_f(r))),
	\end{equation*}
	and this implies $\delta(0,f)=0$, which contradicts $ \delta(0,f)>1/2. $
	Therefore, $ {L}_k[D]g $ must be transcendental. Set \begin{equation}\label{e-44.55}
		\Psi:=\frac{{L}_k[D]g-1}{{L}_k[D]g-1}.
	\end{equation}
	Clearly, $ \Psi $ is a meromorphic function having no zeros and poles, since $ {L}_k[D]f=1 \Leftrightarrow  {L}_k[D]g=1 $ and $ {L}_k[D]=\infty \Leftrightarrow {L}_k[D]g=\infty$. Suppose that $ f_1={L}_k[D]f,\; f_2=\Psi $ and $ f_3=-\Psi {L}_k[D] $. Now \eqref{e-44.55} can be expressed as
	\begin{equation}\label{e-33.77}
		f_1+f_2+f_3\equiv 1.
	\end{equation}
	We also see that 
	\begin{equation}\label{e-44.66}
		\sum_{j=1}^{3}N_{f_j}(r,0)=N_{{L}_k[D]f}(r,0)+N_{{L}_k[D]g}(r,0).
	\end{equation} 
	\noindent By Lemma \ref{lem2.8}, it follows that
	\begin{equation}\label{e-44.77}
		N_{{L}_k[D]f}(r,0)\leq N_f(r,0)+T_{{L}_k[D]f}(r)-T_f(r)+O(\log(rT_f(r))),
	\end{equation}
	holds for all $ r\not\in E. $
	By Lemma \ref{lem2.8}, Lemma \ref{lem2.5} (i), it follows from \eqref{e2.1} and \eqref{e-44.22} that 
	\begin{align}\label{e-44.88}
		N_{{L}_k[D]g}(r,0)&\leq N_g(r,0)+T_{N_{ {L}_k[D]g}}(r)-T_g(r)+O(\log(tT_g(r)))\\&\nonumber\leq N_f(r,0)+N_{{ {L}_k[D]g}}(r,\infty)-N_g(r,0)+O(\log(rT_f(r)))\\&\nonumber\leq N_f(r,0)+k\overline{N}_{{ {L}_k[D]g}}(r,\infty)+O(\log(rT_f(r)))\\&\nonumber\leq N_f(r,0)+k\overline{N}_{{{L}_k[D]f}}(r,\infty)+O(\log(rT_f(r))).
	\end{align} Using \eqref{e-44.66} and \eqref{e-44.77} in \eqref{e-44.88}, we obtain 
	\begin{align}\label{e-44.99}
		\sum_{j=1}^{3}N_{f_j}(r,0)&\leq 2N_f(r,0)+T_{{L}_k[D]f}(r)-T_f(r)+k\overline{N}_{{L}_k[D]f}(r,\infty)\\&\nonumber\quad\quad+O(\log(rT_f(r))).
	\end{align} 
	Next our aim is to prove the functions $ f_1 $, $ f_2 $ and $ f_3 $ are linearly dependent. If this is not the case, then considering Lemma \ref{lem2.10}, Lemma \ref{lem2.5}(ii), and equation \eqref{e-33.33}, we derive the following:
	\begin{align*}
		T_{{L}_k[D]f}(r)&\leq T_f(r)\leq\sum_{j=1}^{3}N_{f_j}(r,0)+2\sum_{j=1}^{3}\overline{N}_{f_{j}}(r,\infty)+O(\log(rT_{f_j}(r)))\\&\leq 2N_f(r,0)+T_{{L}_k[D]f}(r)-T_f(r)+k\overline{N}_{{L}_k[D]f}(r,\infty)+ 2\overline{N}_{{L}_k[D]f}(r,\infty)\\&\quad\quad+2\overline{N}_{{L}_k[D]g}(r,\infty)+O(\log(rT_f(r)))\\&\leq 2N_f(r,0)+T_{{L}_k[D]f}(r)-T_f(r)+(4k+1)\overline{N}_f(r,\infty)+O(\log(rT_f(r))).
	\end{align*} 
	A simple computation shows that
	\begin{equation*}
		T_f(r)\leq 2N_f(r,0)+(4+k)\overline{N}_f(r,\infty)+O(\log(rT_f(r))),
	\end{equation*}
	holds for all $ r\not\in E. $ Thus, we have
	\begin{align*}
		(2\delta(0,f)+(4+k)\Theta(\infty,f)-(k+5))T_f(r)\leq O(\log(rT_f(r)))=o(T_f(r)),
	\end{align*}
	which contradicts $ 2\delta(0,f)+(4+k)\Theta(\infty,f)>k+5. $\vspace{1.5mm}
	
	\noindent Hence, there exist constants $ c_1 $, $ c_2 $ and $ c_3 $, not all zero such that 
	\begin{equation}\label{e-44.1010}
		c_1f_1+c_2f_2+c_3f_3=0.
	\end{equation}
We proceed to distinguish between the following three distinct cases.\vspace{1.5mm}
	
	\noindent{\bf Case 1.}
	Let both the functions $ f_2 $ and $ f_3 $ be non-constants. Our aim is to get a contradiction in this situation. If $ c_1=0 $, then \eqref{e-44.1010} yields that $ c_2f_2+c_3f_3=0. $ If $ c_2=0 $, then $ c_3\neq 0 $, and hence $ f_3=-\Psi{L}_k[D]g\equiv 0 $, which contradicts the assumption that $ f_3 $ is non-constant. If $ c_3=0 $, then $ c_2\neq 0 $, and hence $ f_2\equiv 0 $ which contradicts $ f_2=\Psi $
	was non-constant. Therefore, we must have $ c_2\neq 0 $ and $ c_3\neq 0 $. Therefore, we see that $ f_3=-c_2/c_3 $ \textit{i.e.,} $ {L}_k[D]g=c_2/c_3\neq 0 $, which contradicts the fact that $ {L}_k[D]g $ is transcendental. Therefore, we have $ c_1\neq 0 $. If $ c_3=0 $, then \eqref{e-44.1010} reduces to $ c_1f_1+c_2f_2=0. $ In this case, it is not hard to show that $ c_1\neq 0 $ and $ c_2\neq 0$. Therefore, we have $ f_1=(c_1/c_2)f_2 $ \textit{i.e.,} we have $ {L}_k[D]f=-(c_2/c_1)\Psi\not\equiv 0. $ Again since, $ \Psi $ is a meromorphic function having no zeros and poles, so $ {L}_{k}[D]f $ must be an entire function having no zeros, contradicts the fact that $ {L}_k[D]f(0)=\sum_{j=1}^{k}a_jD^jf(0)=0. $ Therefore, we see that $ c_3\neq 0. $ If $ c_2=0 $, then from \eqref{e-44.1010} follows that $ c_1f_1+c_3f_3=0, $ Since, it is proved $ c_1\neq 0 $ and $ c_3\neq 0 $, it can be seen that 
	\begin{equation}
		\label{e-44.1111} \left(1-\frac{c_1}{c_3}\right)f_1+f_3=1.
	\end{equation} 
It is easy to see that $ c_1\neq c_3 $ in \eqref{e-44.1111}, otherwise $ f_3\equiv 1 $ which leads to a contradiction. Based on the Second Fundamental Theorem of Nevanlinna on $ \mathbb{C}^n, $ and by employing \ref{lem2.5} (ii), we find that
	\begin{align*}
		T_{{L}_k[D]f}(r)=T_{f_1}(r)&< \overline{N}_{f_1}(r,0)+\overline{N}_{f_1}\left(r,\frac{1}{1-\frac{c_1}{c_3}}\right)+\overline{N}_{f_1}(r,\infty)+O(\log(rT_f(r)))\\&\leq \overline{N}_{{L}_k[D]f}(r,0)+\overline{N}_{f_2}\left({r,0}\right)+\overline{N}_{f_1}(r,\infty)+O(\log(rT_f(r)))\\&\leq \overline{N}_{f}(r,0)+T_{{L}_k[D]f}(r)-T_f(r)+\overline{N}_{f}(r,\infty)+O(\log(rT_f(r)))
	\end{align*}
	which implies 
	\begin{align*}
		T_f(r)\leq N_f(r,0)+\overline{N}(r,\infty)+O(\log(rT_f(r)))
	\end{align*}
	for sufficiently large values of $ r $.\vspace{2mm}
	
	\noindent Thus, we obtain
	\begin{equation*}
		(\delta(0,f)+\Theta(\infty,f)-1)T_f(r)\leq O(\log(rT_f(r))),
	\end{equation*}
	which yields $ \delta(0,f)+\Theta(\infty,f)\leq 1 $. Since the condition $ 2\delta(0,f)+(k+4)\Theta(\infty,f)>k+5 $ implies that $ \delta(0,f)+\Theta(\infty,f)>3/2>1, $ we have a contradiction. Thus, since $ c_2\neq 0 $, it follows from equations \eqref{e-33.77} and \eqref{e-44.1010} that
	 \begin{align*}
		(c_1-c_2)f_1+(c_3-c_2)f_3=-c_2.
	\end{align*}
The cases where $ c_1=c_2 $ or $ c_2=c_3 $ lead to a contradiction, as they would respectively imply that $ f_3 $ or $ f_1 $ is constant. Thus, we have $ c_1\neq c_2 $ and $ c_3\neq c_2.$ From the Second Fundamental Theorem for meromorphic functions on $ \mathbb{C}^n $ and Lemma \ref{lem2.5} (ii), we obtain 
\begin{align*}
T_{{L}_k[D]f}(r)=T_{f_1}(r)&< \overline{N}_{f_1}(r,0)+\overline{N}_{f_1}\left(r,\frac{-c_2}{c_1-c_2}\right)+\overline{N}_{f_1}(r,\infty)+O(\log(rT_f(r)))\\&\leq \overline{N}_{{L}_k[D]f}(r,0)+\overline{N}_{f_3}\left({r,0}\right)+\overline{N}_{f}(r,\infty)+O(\log(rT_f(r)))\\&\leq \overline{N}_{{L}_k[D]f}(r,0)+\overline{N}_{{L}_k[D]g}\left({r,0}\right)+\overline{N}_{f}(r,\infty)+O(\log(rT_f(r)))\\&\leq\overline{N}_{f}(r,0)+T_{{L}_k[D]f}(r)-T_f(r)+\overline{N}_{f}(r,0)+k\overline{N}_{f}(r,\infty)\\&\quad+\overline{N}_{f}(r,\infty)+O(\log(rT_f(r)))
\end{align*}
which implies that
\begin{equation*}
T_f(r)\leq 2N_f(r,0)+(k+1)\overline{N}_{f}(r,\infty)+O(\log(rT_f(r))).
\end{equation*}
This further shows that the inequality
	\begin{equation*}
		(2\delta(0,f)+(k+1)\Theta(\infty,f)-(k+2))T_f(r)\leq O(\log(rT_f(r)))
	\end{equation*} 
	implies that $ 2\delta(0,f)+(k+1)\Theta(\infty,f)\leq k+2 $.\vspace{2mm} 
	
	We also note that the condition
	\begin{align*}
		2\delta(0,f)+(k+4)\Theta(\infty,f)>k+5
	\end{align*} 
	shows that
	\begin{align*}
		2\delta(0,f)+(k+1)\Theta(\infty,f)>k+3>k+2,
	\end{align*} 
	which leads to a contradiction. \vspace{1.5mm}
	
	\noindent{\bf Case 2.}
	Let $ f_3=-\Psi{L}_k[D]g $ be a constant. Since $ \Psi $ is a meromorphic function having no zeros and poles, hence $ {L}_k[D]g=-f_3/\Psi $ is an entire function. Thus, we have $ {L}_k[D]g(0)=\sum_{j=1}^{k}a_jD^jg(0)=0, $ and then we see that $ f_3=0, $ and hence $ {L}_k[D]g\equiv 0 $, which contradicts that $ {L}_k[D]g $ is transcendental.
	\vspace{1.2mm}
	
	\noindent{\bf Case 3.} Let $ f_2=\Psi $ be a constant. We recall here $ \Psi $ is a meromorphic function having no zeros and poles, therefore in view of \eqref{e-2.5}, we have
	\begin{equation*}
		\frac{{L}_k[D]f-1}{{L}_k[D]g-1}=c,
	\end{equation*}
	where $ c $ is a non-zero constant. This completes the proof.
\end{proof}
\begin{proof}[\bf Proof of Theorem \ref{th-3.3}]
	We first prove that $ f $ is a polynomial. For the sake of contradiction, assume that $ f $ is a transcendental entire function. Therefore, the function
$ 	{F}(z):={(f(z)-a)}/{b} $ is also transcendental entire. By Lemma \ref{lem2.7}, we obtain
	\begin{equation*}
		T_{{F}(r)}\leq N_{{F}}(r,0)+N_{\mathcal{L}_k[D]F}(r,1)+O(\log(rT_f(r))).
	\end{equation*}
	For positive integers $ j $, it is easy to see that $ D^j{F}=D^jf/b,\; $ $ T_{{F}}(r)=T_f(r)+O(1) $ and $ {L}_{k}[D[f]](r,b)={L}_k[D]f/b. $ thus it follows that
	\begin{equation}\label{e-4.1}
		T_f(r)\leq N_f(r,a)+N_{{L}_k[D]f}(r,b)+O(\log(rT_f(r))).
	\end{equation}
	 The assumption that $ f\neq a $ and $ {L}_k[D]f\neq b $ implies that $ N_f(r,a)=0 $ and $ N_{{L}_k[D]f}(r,b)=0. $ Therefore, equation \eqref{e-4.1} yields 
	\begin{align*}
		T_f(r)=O(\log(rT_f(r))), 
	\end{align*} 
	which contradicts the assumption that $ f $ is transcendental. Thus, $ f $ must be a polynomial. Since $ f(z)\neq a $ for all $ z\in\mathbb{C}, $ it follows from Liouville's theorem that $ f $ must be a constant, which completes the proof.
\end{proof}
\begin{proof}[\bf Proof of Theorem \ref{th-3.5}]
	We first prove that $ f $ is a polynomial. On contrary, we suppose that $ f $ is a transcendental entire function. Set
	\begin{equation*}
		{G}(z):=\frac{1}{b}\left(\beta_m\frac{f^{k+m+1}}{k+m+1}+\beta_{m-1}\frac{f^{k+m}}{k+m}+\cdots+\beta_1\frac{f^{k+1}}{k+1}\right).
	\end{equation*} 
	It is easy to see that $ G $ is also transcendental entire function.  By a simple computation, using \eqref{e-1.1} we obtain  
	\begin{equation}\label{e-4.2}
		DG=\frac{f^kQ_m(f)Df}{b}.
	\end{equation} It is easy to see that zero multiplicity at each point of $ 0 $-divisor of $ {G} $ is at least $ k+1\geq 3. $ Therefore, it follows that 
	\begin{align*}
		\overline{N}_{{G}}(r,0)\leq \frac{1}{3} N_{{G}}(r,0)+O(\log r).
	\end{align*}
	By the assumption, we have 
	\begin{align*}
		DG=\frac{f^kQ_m(f)Df}{b}\neq 1.
	\end{align*}
	In view of Lemma \ref{lem2.11}, we easily obtain
	\begin{align*}
		T_{{G}}(r)&\vspace{4mm}\leq\vspace{4mm} 2\overline{N}_{{G}}(0)+N_{DG}(r,1)+O(\log(rT_f(r)))\vspace{4mm}\\&\vspace{4mm}\leq\frac{2}{3}N_{{G}}(r,0)+O(\log(rT_f(r)))\vspace{4mm}\\&\vspace{4mm}\leq\frac{2}{3}T_{{G}}(r)+O(\log(rT_f(r))) 
	\end{align*} which implies that
	\begin{align*}
		\frac{1}{3}T_{{G}}(r)\leq O(\log(rT_f(r)))
	\end{align*}	which is a contradiction, as $ {G} $ is a transcendental function.\vspace{2mm}
	
	\noindent  Thus, $f$ is a polynomial, and so is $f^kQ_m(f)Df$. The condition $f^kQ_m(f)Df \neq b$ implies that this function must be a constant. If $f$ were non-constant, the degree of $f^kQ_m(f)Df$ would be greater than or equal to 1, contradicting the fact that it is a constant. Therefore, $f$ must be a constant, which completes the proof.
\end{proof}

\vspace{1.6mm}

\noindent\textbf{Compliance of Ethical Standards:}\\

\noindent\textbf{Funding:} The authors have not received any funding  for this research work. \\

\noindent\textbf{Authors Contribution:} Both the authors have prepared the text of the manuscript and both the authors have reviewed the manuscript.\\

\noindent\textbf{Conflict of interest.} The authors declare that there is no conflict  of interest regarding the publication of this paper.\\

\noindent\textbf{Data availability statement.}  Data sharing is not applicable to this article as no datasets were generated or analyzed during the current study.

\end{document}